\documentclass[english,12pt]{article}
\usepackage[T1]{fontenc}
\usepackage[utf8]{inputenc}

\usepackage{fullpage}
\usepackage{amsmath,amsthm}

\newtheorem{thm}{Theorem}
\numberwithin{thm}{section}
\newtheorem{lem}{Lemma}
\numberwithin{lem}{section}

\theoremstyle{definition}

\numberwithin{defi}{section}

\numberwithin{exper}{section}
\theoremstyle{exper}

\theoremstyle{definition}
\newtheorem{assumption}{Assumption}
\numberwithin{assumption}{section}

\theoremstyle{definition}

\numberwithin{exam}{section}

\numberwithin{hequa}{section} 

\newtheorem{rem}{Remark}

\usepackage{amssymb}
\usepackage{graphicx}
\usepackage{subfigure}
\usepackage{color}
\usepackage{epsfig}
\usepackage{epstopdf}
\usepackage{hyperref}
\usepackage{algorithm,algcompatible}
\usepackage{authblk}
\usepackage{comment}

\numberwithin{equation}{section}

\graphicspath{ {./figures/} }

\title{How to avoid order reduction in third-order exponential Runge--Kutta methods for problems with non-commutative operators?}

\author[1]{Thi Tam Dang}
\author[2]{Trung Hau Hoang\thanks{\href{mailto:trunghaugg@gmail.com}{Corresponding author \\
			Email addresses:
			 tam.dang@helsinki.fi (Thi Tam Dang),
			 trunghaugg@gmail.com (Trung Hau Hoang)}}}
\affil[1]{ Department of Mathematics and Statistics, University of Helsinki, Finland} 
\affil[2]{ Department of Mathematics, University of Innsbruck, Austria} 

\date{}
\begin{document}

\maketitle

\section*{Abstract}
This paper investigates the performance of a subclass of exponential integrators, specifically explicit exponential Runge--Kutta methods. It is well known that third-order methods can suffer from order reduction when applied to linearized problems involving unbounded and non-commuting operators. In this work, we consider a fourth-stage third-order Runge--Kutta method, which successfully achieves the expected order of accuracy and avoids order reduction, as long as all required order conditions are satisfied. The convergence analysis is carried out under the assumption of higher regularity for the initial data. Numerical experiments are provided to validate the theoretical results.

\section{Introduction}
Exponential integrators have recently gained attention as efficient numerical methods for solving stiff systems of differential equations (see \cite{ 10.1007/s10543-013-0446-0, CROUSEILLES2020109688, DEKA2023101302,HO2010}). These integrators handle stiffness by solving the linear part exactly while integrating the nonlinear part explicitly. In this paper, we are interested in the nonlinear parabolic equation for the unknown \( u(t, x) \):
\begin{equation}\label{eq98}
	\begin{aligned}
		\partial_{t}u - \triangle u = f(\nabla u, u), \qquad u(0) = u_0,
	\end{aligned}
\end{equation}
where \( f \) is a nonlinear function of \( \nabla u \) and \( u \). For cases where \( f \) depends solely on \( u \), the error analysis of exponential Runge--Kutta methods has been thoroughly studied (see \cite{HO2005}). In this work, we examine the application of exponential Runge--Kutta methods to a simplified linear form of \eqref{eq98}:
\begin{equation}\label{4121}
	u^{\prime}(t) + Au(t) = Bu(t), \quad u(0) = u_0,
\end{equation}
where \( A \) is an operator generating an analytic semigroup, and \( B \) is relatively bounded with respect to \( A \). For instance, \( A \) might be a second-order strongly elliptic operator, while \( B \) could be a first-order differential operator. Understanding the linear problem \eqref{4121} serves as a foundation for analyzing the more complex nonlinear problem \eqref{eq98}.

The exact solution of \eqref{4121} is given by \( u(t) = \mathrm{e}^{-t(A-B)}u_0 \). We treat \( A \) exactly while handling \( B \) explicitly. For example, applying the \textit{variation-of-constants} formula and integrating the nonlinear term at \( u_0 \) leads to the exponential Euler method
\begin{equation*}
	u_{1} = \mathrm{e}^{-\tau A} u_0 + \tau \varphi_1(-\tau A) Bu_0, \qquad \varphi_1(x) = \frac{\mathrm{e}^x - 1}{x}.
\end{equation*}

This approach is motivated by two key factors. 
Firstly, it allows for the investigation of the convergence order of exponential Runge--Kutta methods when applied to unbounded operators such as  \( B \). Secondly, the approximation of the exponential function, given by \( \mathrm{e}^{-tA}u_0 \), is a more computationally feasible approach than the evaluation of the exponential function \( \mathrm{e}^{-t(A-B)}u_0 \) when applied to the initial value problem.

As shown in the preprint \cite{HAU2024}, third-order exponential Runge--Kutta methods suffer a reduction in order to approximately $2.5$ when they fail to satisfy the strong form of the order condition. This observation highlights the need for the development of a method that preserves the desired order of accuracy. In this paper, we propose a fourth-stage Runge–Kutta method that satisfies all the requisite order conditions and, moreover, imposes higher regularity on the initial data. The proposed method successfully overcomes the order reduction issue commonly observed in standard methods, showing significantly improved accuracy and robustness.

The structure of this paper is as follows. Section \ref{recallEXP} provides a review of exponential Runge--Kutta methods and their full order conditions. The convergence analysis of \eqref{4121} is carried out in Section \ref{Analyticalframework}, where the framework of analytic semigroups in a Banach space \( X \) is introduced. In Section \ref{erroranalysis}, we examine the convergence of explicit exponential Runge--Kutta methods applied to \eqref{4121}, with the key results presented in Theorem \ref{theo4}. Numerical investigations are discussed in Section \ref{numericalchap2}. Finally, Section \ref{concluchap2} concludes the paper.

\section{Exponential Runge--Kutta methods}\label{recallEXP}
Let $\tau$ be the time step size, and define \( t_{n+1} = t_n + \tau \) with \( t_n = n\tau \). We define $u_{n}$ as an approximation of $u(t_{n})$ and $ U_{n i}$ as an approximation of $(t_n + c_i \tau)$. Exponential Runge--Kutta methods when applied to \eqref{4121} have the following form
\begin{subequations}\label{eq20}
	\begin{equation}\label{eq20a}
		U_{n i}  =\mathrm{e}^{-c_i \tau A}u_n +\tau \sum_{j=1}^{i-1} a_{i j}\left(-\tau A\right) B  U_{n j}, \quad 1 \leq i \leq s,
	\end{equation}
	\begin{equation}\label{eq20b}
		u_{n+1} =\mathrm{e}^{-\tau A }u_n +\tau\sum_{i=1}^s b_i\left(-\tau A\right) B U_{n i},
	\end{equation}
\end{subequations}
The coefficients \( a_{ij}(-\tau A) \) and \( b_i(-\tau A) \) are typically linear combinations of the functions \( \varphi_k(-\tau A) \), which are defined as

$$
\varphi_k(z) = \int_0^1 e^{(1-\theta) z} \frac{\theta^{k-1}}{(k-1)!} \, d\theta, \quad k \geq 1.
$$

The function \( \varphi_k(-\tau A) \) are well defined by the lemma \ref{parabolicsmoothing1} below  . The coefficients \( a_{ij}(-\tau A) \) and \( b_i(-\tau A) \) ( jointly referred to as \( \phi(-t A) \)) satisfy the following inequality
\begin{equation}\label{eq19}
	\|\phi(-t A)\| + \left\|t^\eta {A}^\eta \phi(-t A)\right\| \leq C, \quad 0 \leq \eta \leq 1.
\end{equation}
We recall the order conditions for the third-order method (see \cite{HO2005}), which satisfies all  required order conditions, as follows:
\begin{subequations} 
	\begin{equation}\label{eq92a}
 b_1(-\tau A) + b_2(-\tau A) + b_3(-\tau A) + b_4(-\tau A) =\varphi_1(-\tau A), \\
	\end{equation}
	\begin{equation}\label{eq92b}
	c_2 b_2(-\tau A)  + c_3 b_3(-\tau A)  + c_4 b_4(-\tau A)   =\varphi_2(-\tau A),
	\end{equation}
	\begin{equation}\label{eq92c}
 a_{21}(-\tau A)= c_2 \varphi_1\left(-c_2 \tau A\right),
\end{equation}
	\begin{equation}\label{eq92d}
	a_{31}(-\tau A) + a_{32}(-\tau A)= c_3 \varphi_1\left(-c_3 \tau A\right),
\end{equation}
	\begin{equation}\label{eq92e}
	a_{41}(-\tau A) + a_{42}(-\tau A) + a_{43}(-\tau A)= c_4 \varphi_1\left(-c_4 \tau A\right),
\end{equation}
	\begin{equation}\label{eq92f}
 c_2^2 b_2(-\tau A)  + c_3^2 b_3(-\tau A)  + c_4^2 b_4(-\tau A) =2\varphi_3(-\tau A),
\end{equation}
	\begin{equation}\label{eq92g}
b_2(-\tau A) J \psi_{2,2}(-\tau A) +b_3(-\tau A) J \psi_{2,3}(-\tau A) +b_4(-\tau A) J \psi_{2,4}(-\tau A)  =0,
\end{equation}
where \( J \) is an arbitrary operator and, 
$$\psi_{2,2}(-\tau A)=c_2^2 \varphi_2\left(-c_2 \tau A\right), $$
$$\psi_{2,3}(-\tau A)=c_3^2 \varphi_2\left(-c_3 \tau A\right)-c_2 a_{32}(-\tau A), $$
$$\psi_{2,4}(-\tau A)=c_4^2 \varphi_2\left(-c_4 \tau A\right)-c_2 a_{42}(-\tau A) - c_3 a_{43}(-\tau A). $$
In the context of the three-stage method, certain terms emerge that are not susceptible to resolution, which could potentially lead to a reduction in the order of the method (see \cite{HAU2024}). It thus follows that a four-stage method is required in order to eliminate these terms and ensure the preservation of third-order accuracy.

\end{subequations}

 \section{Analytical framework}\label{Analyticalframework} 

This section provides the analytical framework that forms the foundation for the upcoming convergence analysis. Let \( X \) be a Banach space, with \( \mathcal{D}(A) \) representing the domain of \( A \) in \( X \). The assumptions on the operators \( A \) and \( B \) follow those in \cite{henry1981geometric,pazy1983semigroups}.

\begin{assumption}\label{ass1}
	Let $A: \mathcal{D}(A) \rightarrow X$ be sectorial.
\end{assumption}
The operator $-A$ is an infinitesimal generator of an analytic semigroup $\left\{\mathrm{e}^{-t A}\right\}_{t \geq 0}$ under Assumption \ref{ass1}.
The stability bounds provided below play a crucial role in the subsequent error analysis (see \cite{HOCHBRUCK2005323}).

\begin{lem}\label{lem1} 
	For fixed $\omega \in (-a,\infty)$ and together with  Assumption \ref{ass1}, the following bounds 
	\begin{equation}\label{parabolicsmoothing1}
		\left\|\mathrm{e}^{-t A}\right\|+\left\|t^\gamma A^\gamma \mathrm{e}^{-t A}\right\| \leq C, \quad \gamma \geq 0, 
	\end{equation}
	hold uniformly on $t \in [0,T]$.
\end{lem}

\begin{assumption}\label{ass2} 
Let \( 0 < \gamma \leq 1 \), and let \( B: \mathcal{D}(B) \rightarrow X \) be a closed linear operator satisfying the following condition:
\begin{equation}\label{eq99}
	\mathcal{D}(A^{\gamma}) \subset \mathcal{D}(B) \quad \text{and} \quad \|B x\| \leq \varepsilon\|A^{\gamma} x\| + K(\varepsilon)\|x\| \quad \text{for all} \ x \in \mathcal{D}(A^{\gamma}),
\end{equation}
for sufficiently small \( \varepsilon > 0 \), where \( K(\varepsilon) \) is a non-negative function. Under Assumption \ref{ass2}, there exists a constant \( C > 0 \) such that (see \cite{pazy1983semigroups})
\begin{equation}\label{eq14}
	\left\| B A^{-\gamma} \right\| \leq C.
\end{equation}

Furthermore, we assume the following condition
\begin{equation}\label{eq15}
	\left\| A^{-\gamma} B \right\| \leq C.
\end{equation}
For the upcoming error analysis, we recall the following lemma from \cite{HAU2024}.
\begin{lem}\label{lem3}
	Under Assumptions \ref{ass1} and \ref{ass2}, we have
	\begin{subequations}
		\begin{equation}\label{eq68a}
			\left\|   \phi \left(-\tau A\right) B\right\|\leq C\tau^{-\gamma},
		\end{equation} 					
		\begin{equation}\label{eq68b}
			\left\|  \mathrm{e}^{-t A} \phi \left(-\tau A\right) B\right\|\leq Ct^{-\gamma},
		\end{equation}
	\end{subequations}
	where $\phi(-\tau A)$ is an arbitrary  linear combinations of the functions $\varphi_k(-\tau A)$.
\end{lem}

\end{assumption}

\section{Convergence results for exponential Runge--Kutta methods}\label{erroranalysis}

In this section, we analyze the convergence analysis of the exponential Runge--Kutta method applied to  \eqref{4121}.

\subsection{Error recursion and defects representation}\label{errorrecursion}
The exact solution of \eqref{4121} can be represented using the variation-of-constants formula as:
\begin{equation}\label{eq91}
	u\left(t_n+\theta \tau \right)=\mathrm{e}^{-\theta \tau A} u\left(t_n\right)+\int_0^{\theta \tau} \mathrm{e}^{-(\theta \tau-\xi) A} Bu\left(t_n+\xi\right) \mathrm{d} \xi .
\end{equation}
The Taylor series expansion of \( u \)
 is expressed as follows:
\begin{equation}\label{eq90} 
	u\left(t_n+\tau\right)=\sum_{j=1}^q \frac{\tau^{j-1}}{(j-1) !} u^{(j-1)}\left(t_n\right)+\int_0^\tau \frac{(\tau-\sigma)^{q-1}}{(q-1) !} u^{(q)}\left(t_n+\sigma\right) \mathrm{d} \sigma .
\end{equation} 
Inserting \eqref{eq90} into the right-hand side of \eqref{eq91}, we obtain
\begin{equation}\label{eq87}
	\begin{aligned}
		u\left(t_n+c_i \tau \right)= & \mathrm{e}^{-c_i \tau A} u\left(t_n\right)+\sum_{j=1}^{q_i}\left(c_i \tau \right)^j \varphi_j\left(-c_i \tau A\right) Bu^{(j-1)}\left(t_n\right) \\
		& +\int_0^{c_i \tau} \mathrm{e}^{-\left(c_i \tau-\xi\right) A} \int_0^\xi \frac{(\tau-\sigma)^{q_i-1}}{\left(q_i-1\right) !} Bu^{\left(q_i\right)}\left(t_n+\sigma\right) \mathrm{d} \sigma \mathrm{d} \xi .
	\end{aligned}
\end{equation}
Substituting the exact solution into the numerical scheme \eqref{eq20} results in
	\begin{subequations}\label{eq84}
	\begin{equation}\label{eq84a}
u\left(t_n+c_i \tau \right)=\mathrm{e}^{-c_i \tau A} u\left(t_n\right)+\tau \sum_{j=1}^{i-1} a_{i j}(-\tau A) Bu\left(t_n+c_j \tau \right)+\Delta_{n i},
	\end{equation}
	\begin{equation}\label{eq84b}
u\left(t_{n+1}\right)=\mathrm{e}^{-\tau A} u\left(t_n\right)+\tau \sum_{i=1}^s b_i(-\tau A) Bu\left(t_n+c_i \tau\right)+\delta_{n+1}
	\end{equation}
\end{subequations}
with defects $\Delta_{n i}$ and $\delta_{n+1}$. Placing \eqref{eq90} in \eqref{eq84a}, we get
\begin{equation}\label{eq86}
	\begin{aligned}
		u\left(t_n+c_i \tau \right)= & \mathrm{e}^{-c_i \tau A} u\left(t_n\right)+ \tau \sum_{k=1}^{i-1} a_{i k}(-\tau A) \sum_{j=1}^{q_i} \frac{\left(c_k \tau \right)^{j-1}}{(j-1) !} Bu^{(j-1)}\left(t_n\right) \\
		& +\tau \sum_{k=1}^{i-1} a_{i k}(-h A) \int_0^{c_k \tau} \frac{\left(c_k \tau-\sigma\right)^{q_i-1}}{\left(q_i-1\right) !} Bu^{\left(q_i\right)}\left(t_n+\sigma\right) \mathrm{d} \sigma+\Delta_{n i} .
	\end{aligned}
\end{equation}

The defects \( \Delta_{n i} \) and \( \delta_{n+1} \) can be explicitly expressed by subtracting \eqref{eq87} from \eqref{eq86}. 
Define the errors between the numerical and exact solutions as \( e_n = u_n - u\left(t_n\right) \) and \( E_{n i} = U_{n i} - u\left(t_n + c_i \tau\right) \). Subtracting \eqref{eq84} from the numerical solution \eqref{eq20} results in the following error recursion
	\begin{subequations}\label{eq85}
	\begin{equation}\label{eq85a}
\quad E_{n i}=\mathrm{e}^{-c_i \tau A} e_n+\tau \sum_{j=1}^{i-1} a_{i j}(-\tau A)B E_{n j}-\Delta_{n i},
	\end{equation}
	\begin{equation}\label{eq85b}
e_{n+1}=\mathrm{e}^{-\tau A} e_n+\tau \sum_{i=1}^s b_i(-\tau A)B E_{n i}-\delta_{n+1}.
	\end{equation}
\end{subequations}
The explicit forms of these error terms will be presented in subsequent sections.

\subsection{Convergence results for third-order methods}
We are now ready to conduct the error analysis of third-order exponential Runge--Kutta methods with
 \( s = 4 \) in \eqref{eq20}. Throughout this section, we assume that \( u_0 \in \mathcal{D}(A^2) \).

\begin{rem}
Using a Taylor series expansion of \( u \), we can differentiate between two cases: \( t_n = 0 \) and \( t_n \neq 0 \), as discussed in \cite{HAU2024}. We have
	\begin{equation}\label{errf32}
		u\left( \xi\right) = u\left( 0 \right)+ \xi u^{\prime} \left( 0 \right)  + \int_0^\xi (\xi - \sigma) u^{ \prime \prime}\left(\sigma\right) \mathrm{d} \sigma ,
	\end{equation}
	and 
	\begin{equation}\label{eq59} 
		u\left( t_n+\xi\right) = u\left( t_n \right)+ \xi u^{\prime} \left( t_n \right)  + \frac{\xi^2}{2} u^{ \prime \prime}(t_n) + \int_0^\xi \frac{(\xi - \sigma)^2}{2} u^{ \prime \prime \prime}\left(t_n+\sigma\right) \mathrm{d} \sigma .
	\end{equation}
The second- and third-order derivatives of $u$ are given by
	\begin{equation}\label{eq16}
		u^{\prime \prime }(t) =  \mathrm{e}^{-t(A-B)} (A-B)^{2}u_0,\, \, \, \,  u^{\prime \prime \prime}(t) =(A-B) \mathrm{e}^{-t(A-B)} (A-B)^{2}u_0.
	\end{equation}
\end{rem}
We use the following lemma to bound $	u^{\prime \prime }(t)$ and $ u^{\prime \prime \prime}(t)$.
\begin{lem}\label{lem5} 
	Under Assumption \ref{ass1} and \( u_0 \in \mathcal{D}(A^{2}) \), we have 
	\begin{subequations}
		\begin{equation}\label{eq73}
			\left\| u^{\prime \prime }(t) \right\|  \leq C , 
		\end{equation} 					
		\begin{equation}\label{eq74}
			\left\| u^{\prime \prime \prime}(t) \right\| \leq C t^{-1}. 
		\end{equation}
	\end{subequations}			 		
\end{lem} 
\begin{proof}
	These bounds can be directly obtained by applying 
	  \eqref{parabolicsmoothing1} and \( u_0 \in \mathcal{D}(A^{2}) \). 
\end{proof}

Following the procedure in \ref{errorrecursion} and using the order conditions \eqref{eq92a}-\eqref{eq92g},  the global error is expressed as
\begin{equation}\label{eq62}
	\begin{aligned}
		e_{n+1} 
		& = \mathrm{e}^{-\tau A} e_n+\tau \mathcal{N} e_n+\tau^2 \mathcal{M} e_n+\tau^3 \mathcal{K} e_n+\tau^4 \mathcal{P} e_n+ \mathcal{T}_{n+1} + \mathcal{R}_{n+1} + \mathcal{ \widehat{R}}_{n+1} , 
	\end{aligned}
\end{equation}
where the terms independent of $n$ are provided by
\begin{equation*}
	\begin{aligned}
		\mathcal{N} &=  \sum_{q=1}^4   b_i(-\tau A) B \mathrm{e}^{-c_i \tau A},    \\
		\mathcal{M}  & = b_2(-\tau A) B      a_{21}\left(-\tau A\right)  B + \sum_{q=1}^2 b_3(-\tau A) B a_{3q}\left(-\tau A\right) B \mathrm{e}^{-c_q \tau A} \\
		& \qquad + \sum_{q=1}^4 b_4(-\tau A) B a_{4q}\left(-\tau A\right) B \mathrm{e}^{-c_q \tau A} , \\
		\mathcal{K}  & =  \sum_{q=3}^4 b_q(-\tau A) B a_{q2}\left(-\tau A\right) B a_{21}\left(- \tau A\right)  B e_n + \sum_{q=1}^2 b_4(-\tau A)B a_{43} \left( -\tau A\right)B a_{3q} \left( -\tau A\right)Be_n, \\
		\mathcal{P}  & =   b_3(-\tau A) B       a_{43} \left( -\tau A\right)B a_{32}\left( -\tau A\right) B a_{21}\left(- \tau A\right)  Be_n. 
	\end{aligned}
\end{equation*}
Additionally, the terms that depend on $n$ are represented by
\begin{equation*}
	\begin{aligned}
		\mathcal{T}_{n+1} & =  -\tau^4 b_3(-\tau A) B a_{32}\left( -\tau A \right) B \psi_{2,2}(-\tau A) B u^{\prime}(t_n) \\ & \qquad - \tau^4 \sum_{j=2}^{3} b_4(-\tau A) B a_{4j} \left( -\tau A \right) B \psi_{2,j}(-\tau A) B u^{\prime}(t_n)  \\
		& \qquad - \tau^5 b_4(-\tau A) B a_{43} \left( -\tau A\right)B a_{32}\left( -\tau A\right)B \psi_{2,2}(-\tau A) Bu^{\prime}\left(t_n\right),  \\
		\mathcal{R}_{n+1} 	& = \sum_{j=0}^{n-1}    \mathcal{R}_{n+1}^{[j]} , \qquad 1\leq j \leq 11, \\
		\mathcal{\widehat{R}}_{n+1} 	& = \sum_{j=0}^{n-1}    \mathcal{\widehat{R}}_{n+1}^{[j]} , \qquad 1\leq j \leq 4, 		 
	\end{aligned}
\end{equation*}
where
{\allowdisplaybreaks
	\begin{align*}
		\mathcal{R}_{n}^{[k]}	& =  - \tau b_{k+1}(-\tau A)B   \int_0^{c_{k+1} \tau} \mathrm{e}^{-(c_{k+1} \tau-\xi) A} \int_0^\xi (\xi - \sigma) Bu^{\prime \prime}\left(t_{n-1}+\sigma\right) \mathrm{d} \sigma   \mathrm{d} \xi , \quad  k =1,2,3,
		\\
		\mathcal{R}_{n}^{[4]} & = - \tau^2 b_3(-\tau A) B    a_{32}\left( -\tau A\right) B   \int_0^{c_2 \tau} \mathrm{e}^{-(c_2 \tau-\xi) A} \int_0^\xi (\xi - \sigma) Bu^{\prime \prime}\left(t_{n-1}+\sigma\right) \mathrm{d} \sigma   \mathrm{d} \xi ,    \\
		\mathcal{R}_{n}^{[l+3]} & = - \tau^2 b_4(-\tau A)B  a_{4l}(-\tau A)B \int_0^{c_{l} \tau} \mathrm{e}^{-(c_{l} \tau-\xi) A} \int_0^\xi (\xi - \sigma) Bu^{\prime \prime}\left(t_{n-1}+\sigma\right) \mathrm{d} \sigma   \mathrm{d} \xi ,  \quad l =2,3, \\
		\mathcal{R}_{n}^{[7]} & =  \tau^2 b_3(-\tau A) B  a_{32}\left(-\tau A\right) B \int_0^{c_2\tau} \left(c_2\tau - \sigma\right) Bu^{\prime \prime}\left(t_{n-1}+\sigma\right) \mathrm{d} \sigma  ,   \\
		\mathcal{R}_{n}^{[m+6]} & = - \tau^2 b_4(-\tau A)B  a_{4m}(-\tau A)B   \int_0^\xi (c_{m} \tau - \sigma) Bu^{\prime \prime}\left(t_{n-1}+\sigma\right) \mathrm{d} \sigma  ,  \quad m =2,3, \\	
		\mathcal{R}_{n}^{[10]} & = \tau^3 b_4(-\tau A)B  a_{43}(-\tau A)B  a_{32}\left(-\tau A\right) B \int_0^{c_2\tau} \left(c_2\tau - \sigma \right) Bu^{\prime \prime}\left(t_{n-1}+\sigma\right) \mathrm{d} \sigma,    \\		
		\mathcal{R}_{n}^{[11]} & = - \tau^3 b_4(-\tau A)B  a_{43}(-\tau A)B  a_{32}\left( -\tau A\right) B   \int_0^{c_2 \tau} \mathrm{e}^{-(c_2 \tau-\xi) A} \int_0^\xi (\xi - \sigma) Bu^{\prime \prime}\left(t_{n-1}+\sigma\right) \mathrm{d} \sigma   \mathrm{d} \xi,    \\										
		\mathcal{\widehat{R}}_{i}^{[1]} & = 		\int_0^{ \tau} \mathrm{e}^{-(\tau-\xi) A} \int_0^\xi \frac{ (\xi - \sigma)^2}{2}  Bu^{ \prime \prime \prime}\left(t_{i-1}+\sigma\right) \mathrm{d} \sigma   \mathrm{d} \xi , \qquad n-1 \geq i \geq 2,\\
		\mathcal{\widehat{R}}_{i}^{[k]} & =- \tau b_{k}\left(-\tau A\right) \int_0^{c_{k}\tau} \frac{(c_{k}\tau - \sigma)^2}{2} Bu^{\prime \prime \prime}\left(t_{i-1}+\sigma\right) \mathrm{d} \sigma , \qquad n-1 \geq i \geq 2, \qquad k = 2,3,4,
	\end{align*}
}
and 
\begin{equation*}
	\begin{aligned}
		\mathcal{\widehat{R}}_{1}^{[1]} & = \int_0^{ \tau} \mathrm{e}^{-(\tau-\xi) A} \int_0^\xi  (\xi - \sigma)  Bu^{\prime  \prime}\left(\sigma\right) \mathrm{d} \sigma   \mathrm{d} \xi, \\
		\mathcal{\widehat{R}}_{1}^{[k]} & = - \tau b_k\left(-\tau A\right) \int_0^{c_k\tau} (c_k\tau - \sigma ) Bu^{ \prime \prime}\left(\sigma\right) \mathrm{d} \sigma, \qquad k = 2,3,4.
	\end{aligned}
\end{equation*}  

Solving the error recursion \eqref{eq62} yields
\begin{equation*}
	\begin{aligned}
		e_n& = \tau \sum_{j=0}^{n-1} \mathrm{e}^{-(n-j-1) \tau A} \mathcal{N} e_j + \tau^2 \sum_{j=0}^{n-1} \mathrm{e}^{-(n-j-1) \tau A} \mathcal{M} e_j+ \tau^3 \sum_{j=0}^{n-1} \mathrm{e}^{-(n-j-1) \tau A} \mathcal{K} e_j\\
		& \qquad + \tau^3 \sum_{j=0}^{n-1} \mathrm{e}^{-(n-j-1) \tau A} \mathcal{P} e_j +\tau^3 \sum_{j=0}^{n-1} \mathrm{e}^{-j \tau A} \mathcal{T}_{n-j} + \sum_{j=0}^{n-1} \mathrm{e}^{-j \tau A} \mathcal{R}_{n-j} + \sum_{j=0}^{n-1} \mathrm{e}^{-j \tau A} \mathcal{ \widehat{R}}_{n-j} .
	\end{aligned}
\end{equation*}
The bound of \( e_n \) can be derived by estimating each individual term in the expression for 
 \( e_n \).
\subsection{Preliminary error bounds}\label{error-analysis}

		The bound of $\tau \sum_{j=0}^{n-1} \mathrm{e}^{-(n-j-1) \tau A} \mathcal{N} e_j$ can be obtained using \eqref{eq68a}, \eqref{eq68b}, and \eqref{parabolicsmoothing1} as follows
\begin{equation}\label{er6} 
	\begin{aligned}
		\left \| \tau \sum_{j=0}^{n-1} \mathrm{e}^{-(n-j-1) \tau A} \mathcal{N} e_j  \right \| 
		& \leq    C \tau \sum_{j=0}^{n-2} \left \| \mathrm{e}^{-(n-j-1) \tau A}  \sum_{q=1}^4 b_q(-\tau A)  B   \right\| \left \| e_j \right \|
		\\ 
		& \qquad +  C \tau \left \| \sum_{q=1}^4 b_q(-\tau A)  B \right\|   \left \| e_{n-1}      \right \| \\
		& \leq   C \tau \sum_{j=0}^{n-2} t_{n-j-1}^{-\gamma}   \left \| e_j \right \| +  C \tau^{1-\gamma} \left \| e_{n-1}     \right \|.					   
	\end{aligned}
\end{equation}	
We observe that the estimate of $\tau^2 \sum_{j=0}^{n-1} \mathrm{e}^{-(n-j-1) \tau A} \mathcal{M} e_j$ is almost the same as \\$\tau \sum_{j=0}^{n-1} \mathrm{e}^{-(n-j-1) \tau A} \mathcal{N} e_j$. The only difference is that the order of the former will be an order of $\gamma$ higher than that of the latter due to the additional term $\tau  a_{ij}\left( -\tau A\right) B $ which can be estimated by $$\tau \left \| a_{ij}\left( -\tau A\right) B \right \| \leq \tau^{1-\gamma}.$$ 
This leads to the following estimate 
\begin{equation}
	\begin{aligned}
		\left \| \tau \sum_{j=0}^{n-1} \mathrm{e}^{-(n-j-1) \tau A} \mathcal{M} e_j  \right \| 
		& \leq    \tau^{1-\gamma} \left( C \tau \sum_{j=0}^{n-2} t_{n-j-1}^{-\gamma}   \left \| e_j \right \| +  C \tau^{1-\gamma} \left \| e_{n-1}     \right \| \right) \\
		& \leq   C \tau \sum_{j=0}^{n-2} t_{n-j-1}^{-\gamma}   \left \| e_j \right \| +  C \tau^{1-\gamma} \left \| e_{n-1}     \right \| .					   				
	\end{aligned}
\end{equation}	
Here we assume that the time step size $\tau < 1$, which is a reasonable choice as we are primarily concerned with the asymptotic behavior of the method.
 We notice that the same argument also apply to estimate $\tau^3 \sum_{j=0}^{n-1} \mathrm{e}^{-(n-j-1) \tau A} \mathcal{K} e_j$ and $\tau^3 \sum_{j=0}^{n-1} \mathrm{e}^{-(n-j-1) \tau A} \mathcal{P} e_j$. For brevity, we omit the detailed calculations here and proceed to establish the following auxiliary lemma.

		\begin{lem}\label{lem8}
Under Assumptions \ref{ass1} and \ref{ass2}, and given that \( \| A^{\widetilde{\Gamma}} BA^{-2} \| \leq C \), where \( \widetilde{\Gamma} \in \mathbb{R} \) is fixed and \( u_0 \in \mathcal{D}(A^2) \), the following estimate satisfied:

\[
\| B \phi(-\tau A) B u^{\prime}(t) \| \leq t^{-1} \tau^{\widetilde{\Gamma} - \gamma},
\]

where \( \phi(-\tau A) \) is a linear combination of the functions \( \varphi_k(-\tau A) \).
\end{lem}
\begin{proof}
	 The bound can be derived as follows:
	\begin{align*}
		\left\|  B \phi\left(-\tau A\right) B u^{\prime}\left(t\right) \right\| 
		& \leq \left \|  B A^{-\gamma } \right \| \left \|\phi\left(-\tau A\right)A^{-\widetilde{\Gamma}+\gamma} \right \| \left \| A^{\widetilde{\Gamma}} BA^{-2} \right \| \left \|A^2 (A-B)^{-2} \right \| \left \| u^{\prime \prime \prime}\left(t\right) \right \| \\
		& \leq t^{-1} \tau^{\widetilde{\Gamma}-\gamma} .
	\end{align*}
Here, we employed \eqref{eq14}, \eqref{eq19}, and Lemma \ref{lem5}.
\end{proof}

We are ready to estimate \( \sum_{j=0}^{n-1} \mathrm{e}^{-j \tau A} \mathcal{T}_{n-j} \). We first observe that the estimate of the last term of \(\mathcal{T}_{n}\) always has a higher order compared to the other terms. Using \eqref{eq68a},  we have
\begin{align*}
	 C \tau^5 & \left \|  b_4(-\tau A) B \right \| \left \| a_{43} \left( -\tau A\right)B  a_{32}  \left( -\tau A\right)   B \psi_{2,2}(-\tau A) Bu^{\prime}\left(t_n\right) \right \| \\	\leq 
	 &  C \tau^{5-\gamma}  \left \| a_{43} \left( -\tau A\right)B a_{32}\left( -\tau A\right)   B \psi_{2,2}(-\tau A) Bu^{\prime}\left(t_n\right) \right \|.
\end{align*}
Therefore, the last term of $\mathcal{T}_{n}$ can be omitted. Employing Lemma \ref{lem8} and \eqref{eq68a}, the term for $j = 0$ can be bounded by
{\allowdisplaybreaks
	\begin{align*}
		\left \| \mathcal{T}_{n}  \right \|  
		& \leq  C \tau^4 \left \| b_3(-\tau A) B \right \| \left \|  a_{32}\left( -\tau A \right) \right \| \left \| B  \psi_{2,2}(-\tau A) B  u^{\prime}\left(t_{n-1}\right) \right \| \\
		& \qquad + C \tau^4  \sum_{j=2}^{3} \left \| b_4(-\tau A) B \right \| \left \| a_{4j}\left( -\tau A \right) \right \| \left \|  B \psi_{2,j}(-\tau A) B u^{\prime}(t_{n-1}) \right \| \\
		&  \leq C t_{n-1}^{-1} \tau^{4- 2\gamma + \widetilde{\Gamma}} .
	\end{align*}
}In addition to \eqref{eq68b}, the term for $j = n-1$ can be bounded by	
{\allowdisplaybreaks
	\begin{align*}
		\left \| \mathrm{e}^{-(n-1) \tau A}  \mathcal{T}_{1} \right \|  
		& \leq  C \tau^4  \left \| \mathrm{e}^{-(n-1) \tau A}  b_3(-\tau A) B \right \|  \left \| a_{32}\left( -\tau A \right) B  \right \| \left \|  \psi_{2,2}(-\tau A) \right \| \left \|B  u^{\prime}\left(0\right) \right \|\\
		& \qquad + C  \tau^4 \sum_{j=2}^{3} \left \| \mathrm{e}^{-(n-1) \tau A}  b_4(-\tau A) B \right \| \left \| a_{4j}\left( -\tau A \right) B \right \|  \left \|  \psi_{2,j}(-\tau A) \right \| \left \|B  u^{\prime}\left(0\right) \right \|  \\	
		&  \leq C t_{n-1}^{-\gamma} \tau^{4 - \gamma+ \widetilde{\Gamma}}  .		
	\end{align*}
}The sum of the remaining terms with $j \neq 0$ and $j \neq n-1$ can be written as follows:
\begin{equation}\label{eq21}
	\begin{aligned}
		& \text{I} = - \tau^{4} \sum_{j=1}^{n-2}  \left( \mathrm{e}^{-j \tau A}   A^{\gamma } \right)   \left(b_3(-\tau A) \right) \left( A^{-\gamma } B  \right) \left( a_{32}(-\tau A) \right) \left ( B  \psi_{2,2} \left(- \tau A\right) B u^{\prime }\left(t_{n-j-1}\right) \right )  \\
		& \quad	-\sum_{l=2}^{3}  \tau^{4} \sum_{j=1}^{n-2}  \left( \mathrm{e}^{-j \tau A}   A^{\gamma } \right)  \left(b_4(-\tau A) \right) \left( A^{-\gamma } B  \right) \left( a_{4j}(-\tau A) \right)  \left ( B \psi_{2,l} (-\tau A)  B u^{\prime }\left(t_{n-j-1}\right) \right ) .
		\end{aligned} 	
\end{equation}
By applying \eqref{parabolicsmoothing1} and Lemma \ref{lem8} to estimate each term inside the parentheses in \eqref{eq21}, we derive the following result
{\allowdisplaybreaks
	\begin{align*}
		\left \| \text{I} \right \|
		& \leq C \tau^{ 4 - \gamma + \widetilde{\Gamma}}  \sum_{j=1}^{n-2}	t_{j}^{- \gamma } t_{n-j-1}^{-1} \leq C  \tau^{3 - \gamma + \widetilde{\Gamma} - \zeta} .	\end{align*}		
}By setting $\widetilde{\Gamma} =   \frac{1}{2}- \zeta$, we obtain the order of convergence as $ \frac{7}{2} -\gamma - 2 \zeta$ for the estimation of $ \sum_{j=0}^{n-1} \mathrm{e}^{-j \tau A} \mathcal{T}_{n-j}$. Throughout the paper, \( \zeta \)   denotes a fixed, arbitrarily small positive constant.

We now proceed to bound \( \sum_{j=0}^{n-1} \mathrm{e}^{-j \tau A} \mathcal{R}_{n-j} \).  The estimates for 
\( \mathcal{R}_{n}^{[1]} \), \( \mathcal{R}_{n}^{[2]} \), and \( \mathcal{R}_{n}^{[3]} \) as well as for \( \mathcal{R}_{n}^{[4]} \), \( \mathcal{R}_{n}^{[5]} \), and \( \mathcal{R}_{n}^{[6]} \), follow a similar approach. Thus, we focus on estimating a single term. By applying \eqref{eq68a} and \eqref{eq73}, we can derive bounds for \( \mathcal{R}_{n}^{[1]} \) and \( \mathcal{R}_{n}^{[4]} \) when \( j=0 \) as follows:
{\allowdisplaybreaks
	\begin{align*}
		\left \|   \mathcal{R}_{n}^{[1]}  \right \| 
		& \leq  \tau \left \| b_2(-\tau A) B  \right \|  \int_0^{c_2 \tau} \int_0^\xi \left \| \mathrm{e}^{-\left( c_2 \tau-\xi \right) A} B    \right \| (\xi - \sigma)  \left \| u^{ \prime \prime}\left(t_{n-1}+\sigma\right) \right \| \mathrm{d} \sigma   \mathrm{d} \xi  \\								
		& \leq 	C  \tau^{1 - \gamma} 	\int_0^{c_2 \tau} \int_0^\xi  \left( c_2 \tau-\xi \right)^{- \gamma} (\xi - \sigma)  \mathrm{d} \sigma   \mathrm{d} \xi  \leq 	C  \tau^{4- 2\gamma} ,\\
		\left \| \mathcal{R}_{n}^{[4]} \right \|
		& \leq C \tau^2 \left\|  b_3(-\tau A)  B  \right \| \left\|  a_{32}\left( -\tau A\right) B \right \| \times \\ & \qquad \int_0^{c_2 \tau} \int_0^\xi  \left \| \mathrm{e}^{-\left( c_2 \tau-\xi \right) A} B    \right \|   (\xi - \sigma)  \left \| u^{ \prime \prime}\left(t_{n-1}+\sigma\right)  \right \|  \mathrm{d} \sigma   \mathrm{d} \xi  \\				  
		& \leq 	C  \tau^{2-2\gamma} 	\int_0^{c_2 \tau} \int_0^\xi \left( c_2 \tau-\xi \right)^{-\gamma} (\xi - \sigma)  \mathrm{d} \sigma   \mathrm{d} \xi  \leq 	C \tau^{5-3\gamma}	.	 		
	\end{align*}
} We note that the estimate for \( \mathcal{R}_{n}^{[4]} \) is very similar to that for \( \mathcal{R}_{n}^{[1]} \), with the main difference is that the order of \( \mathcal{R}_{n}^{[4]} \) is $1-\gamma$ higher order  than that of \( \mathcal{R}_{n}^{[1]} \) due to the additional term \( \tau \left\| a_{32}\left(-\tau A\right) B \right\| \leq \tau^{1-\gamma} \). Therefore, we restrict our attention to the estimation of terms containing \( \mathcal{R}_{l}^{[1]} \), where \( 1 \leq l \leq n-1 \). In a similar, we proceed further with the estimation of the remainnng term
\begin{equation*}
	\begin{aligned}
		\left \| \mathcal{R}_{n}^{[7]}  \right \| 
		& \leq \tau^2    \left \| b_3(-\tau A) B \right \| \left\|  a_{32}\left(-\tau A\right) B \right \| \int_0^{c_2\tau} \left(c_2\tau - \sigma \right)  \left\| u^{\prime \prime}\left(t_{n-1}+\sigma\right)  \right \| \mathrm{d} \sigma     \\				
		& \leq 	C  \tau^{2-2\gamma} 	\int_0^{c_2\tau} \left(c_2\tau - \sigma \right)  \mathrm{d} \sigma  \leq 	C  \tau^{4-2\gamma} 	.	
	\end{aligned}
\end{equation*}
We observe that the estimates for \( \mathcal{R}_{n}^{[1]} \) and \( \mathcal{R}_{n}^{[7]} \) are identical, allowing us to perform the estimation only once. This similarity arises because, in \( \mathcal{R}_{n}^{[1]} \), the additional integral, combined with \( \left\| \mathrm{e}^{-\left( c_2 \tau - \xi \right) A} B \right\| \) leads to an estimate of the form \( 1 - \gamma \). For \( \mathcal{R}_{n}^{[7]} \), the additional term \( \tau a_{32}\left(-\tau A\right) B \) also results in the same estimate of \( 1 - \gamma \). Furthermore, the estimates for \( \mathcal{R}_{n}^{[7]} \), \( \mathcal{R}_{n}^{[8]} \), and \( \mathcal{R}_{n}^{[9]} \) can be derived in the same way. Based on these observations, the estimation of \( \mathcal{R}_{n}^{[10]} \) and \( \mathcal{R}_{n}^{[11]} \) can be omitted, allowing us to focus solely on the term
 \( \mathcal{R}_{n}^{[1]} \).

To estimate \( \mathcal{R}_{n}^{[1]} \) for $j=n-1$, we use \eqref{eq68a}, \eqref{eq68b}, and \eqref{eq73} yielding the following:
\begin{equation*}
	\begin{aligned}
		\left \|\mathrm{e}^{-(n-1) \tau A}   \mathcal{R}_{1}^{[1]}  \right \| 
		& \leq \tau \left \| \mathrm{e}^{-(n-1) \tau A}  b_2(-\tau A)  B \right \|     \int_0^{c_2 \tau} \int_0^\xi  \left \|  \mathrm{e}^{-\left(c_2 \tau-\xi \right) A} B \right\|  (\xi - \sigma)   \left \| u^{ \prime \prime}\left(\sigma\right) \right \| \mathrm{d} \sigma   \mathrm{d} \xi  \\									
		& \leq 	C t_{n-1}^{ - \gamma} \tau 	\int_0^{c_2 \tau} \int_0^\xi \left( c_2 \tau-\xi \right)^{-\gamma}(\xi - \sigma)  \mathrm{d} \sigma   \mathrm{d} \xi   \leq 	C  t_{n-1}^{ - \gamma} \tau^{4-\gamma} .		
	\end{aligned}
\end{equation*}
After handling the boundary terms separately, the remaining terms in the sum, where \( j \neq 0 \) and \( j \neq n-1 \), can be bounded using \eqref{parabolicsmoothing1}, \eqref{eq68a}, \eqref{eq68b}, and \eqref{eq74}, as shown below:
\begin{equation*}
	\begin{aligned}
		\left \| \sum_{j=1}^{n-2} \mathrm{e}^{-j \tau A}  \mathcal{R}_{n-j}^{[1]}  \right \| 
		& \leq  \tau  \sum_{j=1}^{n-2} \left\| \mathrm{e}^{-j \tau A}  b_2(-\tau A)  B \right\| \times   \\
		& \qquad \int_0^{c_2 \tau}\int_0^\xi  \left \| \mathrm{e}^{-\left( c_2 \tau-\xi \right) A} \right\| \left \| A^{-\beta}  \right\| \left \| A^{\beta} B A^{-1} \right\|  (\xi - \sigma)  \left \| u^{ \prime \prime \prime}\left(t_{n-j-1}+\sigma\right) \right \|  \mathrm{d} \sigma   \mathrm{d} \xi  \\					
		& \leq 	C \tau  \sum_{j=1}^{n-2} t_{j}^{-\gamma} t_{n-j-1}^{-1} 	\int_0^{c_2 \tau} \int_0^\xi  (\xi - \sigma)  \mathrm{d} \sigma   \mathrm{d} \xi  \leq 	C t_n^{-\gamma +\zeta} \tau^{3-\zeta}.
	\end{aligned}
\end{equation*}
Here, \( \beta \in \mathbb{R} \) is a fixed constant, and it is also assumed that\( \left\| A^\beta B A^{-1} \right\| \) is bounded.

Lastly, we bound \( \sum_{j=0}^{n-1} \mathrm{e}^{-j \tau A} \mathcal{\widehat{R}}_{n-j} \). Based on the discussion regarding the estimate of \( \sum_{j=0}^{n-1} \mathrm{e}^{-j \tau A} \mathcal{R}_{n-j} \),  it is sufficient to focus on \( \mathcal{\widehat{R}}_{n}^{[1]} \) for \( n = \overline{1, n-1} \). The estimates of \( \mathcal{\widehat{R}}_{n}^{[k]} \), where \( k = 1, 2, 3 \) and \( n = \overline{1, n-1} \), follow the same approach as those for \( \mathcal{\widehat{R}}_{n}^{[1]} \) for \( n = \overline{1, n-1} \). For $j=0$, using \eqref{eq68a} and \eqref{eq74}, we get 
{\allowdisplaybreaks
	\begin{align*}
		\left \| \mathcal{\widehat{R}}_{n}^{[1]}  \right \|
		& \leq    \int_0^{ \tau} \int_0^\xi   \left \| \mathrm{e}^{-\left(  \tau-\xi \right) A} B   \right \|  \frac{ (\xi - \sigma)^2}{2}  \left \| u^{\prime \prime \prime}\left(t_{n-1}+\sigma\right) \right \|  \mathrm{d} \sigma   \mathrm{d} \xi  \\				 
		& \leq 	C  t_{n-1}^{-1}	\int_0^{ \tau} \int_0^\xi  \left(  \tau-\xi \right)^{-\gamma}  \frac{ (\xi - \sigma)^2}{2}  \mathrm{d} \sigma   \mathrm{d} \xi  \leq 	C t_{n-1}^{-1} \tau^{4-\gamma} 	.	
	\end{align*}
}

At \( j = n-1 \), we get \( t_n = 0 \). Thus, by applying \eqref{parabolicsmoothing1} and \eqref{eq73}, we obtain
\begin{equation*}
	\begin{aligned}
		&\left \| \mathrm{e}^{-(n-1) \tau A} \mathcal{\widehat{R}}_{1}^{[1]}  \right \|\leq  \left \| \mathrm{e}^{-(n-1) \tau A}  \right \|   \int_0^{ \tau} \int_0^\xi \left\| \mathrm{e}^{-\left(  \tau-\xi \right) A} \right\|    (\xi - \sigma) \left \| u^{ \prime \prime}\left(\sigma\right) \right \|  \mathrm{d} \sigma   \mathrm{d} \xi  \\	&
		\qquad \qquad \qquad \quad \leq 	C  	\int_0^{ \tau} \int_0^\xi  (\xi - \sigma)    \mathrm{d} \sigma   \mathrm{d} \xi  \leq 	C   \tau^{3}.									 		
	\end{aligned}
\end{equation*}
 For \( j \neq 0 \) and \( j \neq n-1 \), the remaining sum can be bounded using \eqref{parabolicsmoothing1} and \eqref{eq74}
{\allowdisplaybreaks
	\begin{align*}
		& \left \| \sum_{j=1}^{n-2} \mathrm{e}^{-j \tau A}  \mathcal{R}_{n-j}^{[5]}   \right \|  \leq  \sum_{j=1}^{n-2} \left\| \mathrm{e}^{-j \tau A} A^{1-\zeta} \right \|  \int_0^{ \tau} \int_0^\xi  \left\| \mathrm{e}^{-\left(  \tau-\xi \right) A} A^{\zeta} \right \|    \frac{ (\xi - \sigma)^2}{2}  \left \|  u^{\prime \prime \prime}\left(t_{n-j-1}+\sigma\right) \right \|  \mathrm{d} \sigma   \mathrm{d} \xi  \\					 
		&  \qquad \qquad \qquad \quad \, \, \, \, \, \leq	C  \sum_{j=1}^{n-2} t_{j}^{\zeta-1} t_{n-j-1}^{-1}	\int_0^{ \tau} \int_0^\xi  \frac{ (\xi - \sigma)^2}{2} \left(  \tau-\xi \right)^{-\zeta} \mathrm{d} \sigma   \mathrm{d} \xi \\
		&   \qquad \qquad \qquad \quad \, \, \, \, \, \leq 	C \tau^{4-2\zeta} \sum_{j=1}^{n-2} t_{j}^{\zeta-1} t_{n-j-1}^{-1 +\zeta}    \leq 	C t_n^{-1+2\zeta } \tau^{3-2\zeta}  .
	\end{align*}
}

 We are ready to state the main result.
\subsection{Error bounds} 

			\begin{thm}\label{theo4} 
	Let Assumptions \ref{ass1} and \ref{ass2} hold.  Furthermore, assume that \( \left\| A^\beta B A^{-1} \right\| \), \( \left\| A^{-2} B A^{\widetilde{\Gamma}} \right\| \), where $\beta, \widetilde{\Gamma} \in \mathbb{R}$, are bounded, and that \( u_0 \in \mathcal{D}(A^{2}) \). Additionally, suppose all the order conditions are satisfied in the strong form. Under these conditions, the numerical solution of the initial value problem \eqref{4121}, obtained using third-order exponential Runge–Kutta methods with $s=4$ in \eqref{eq20}, satisfies the error bound:
	\[
	\left\|u_n - u\left(t_n\right)\right\| \leq C t_n^{-1+2\zeta} \tau^{3-2\zeta}
	\]
	uniformly for $0 \leq n\tau \leq T$. The constant $C$ depends on $T$, but it is independent of $n$ and $\tau$.
\end{thm}
\begin{proof}
	By combining the estimates of each term in \( e_n \) and using the triangle inequality, the bound for \( e_n \) can be expressed as follows:
	\begin{equation*}
		\begin{aligned}
			\left \| e_n \right \| 
			& \leq  C \tau \sum_{j=1}^{n-2} t_{n-j-1}^{-\gamma}   
			\left \| e_j \right \|  + C \tau^{1-\gamma}     \left \| e_{n-1} \right \|  + C t_n^{-1+2\zeta} \tau^{3-2\zeta}  .						
		\end{aligned}
	\end{equation*}
	The proof is completed by applying a discrete Gronwall Lemma (see \cite{HO2010}).
\end{proof}
Throughout the analysis, the operators \( A \) and \( B \) are considered as general operators. In the finite dimensional case, we consider specific operators where \( A \) is a second-order differential operator with homogeneous Dirichlet boundary conditions, and \( B \) is a first-order differential operator. In this scenario, the bounds for \( \left\| A^{-2} B A^{\widetilde{\Gamma}} \right\| \) and \( \left\| A^\beta B A^{-1} \right\| \) are provided by Lemma \ref{Fourierlemma} and Lemma 4.8 in \cite{HAU2024}, which confirm the estimates in Section \ref{error-analysis}.

\begin{lem}\label{Fourierlemma}
	Let  \(\Omega\) be a bounded domain in \(\mathbb{R}^N\), and let \(\zeta > 0\) be a fixed small positive constant. 
	 Consider \( A \) as a second-order differential operator with homogeneous Dirichlet boundary conditions, and \( B \) as a first-order differential operator. In this case, \( A^{-2}BA^{\widetilde{\Gamma}} \) is bounded for the following values of \(\widetilde{\Gamma} \): for \(\widetilde{\Gamma} = 1 - \zeta\) in \(L^1(\Omega)\), \(L^{\infty}(\Omega)\), and for \(\widetilde{\Gamma} = \frac{5}{4} - \zeta\) in \(L^2(\Omega)\).
\end{lem}

\begin{proof}
	First, we observe that, according to the assumptions about the operators \( A \) and \( B \) in this lemma, they are not commutative. Without loss of generality, we write \( f \) using its Fourier series as:
	\begin{equation}\label{eq22}
f(\mathbf{x})=\sum_{\mathbf{k} \in \mathcal{I}_{\mathbf{M}}} \hat{f}(\mathbf{k}) \prod_{i=1}^N \sin \left(k_i \pi x_i\right),
	\end{equation}
where $\mathcal{I}_{\mathbf{M}}=\left\{\mathbf{k}=\left(k_1, k_2, \ldots, k_N\right) \mid 1 \leq k_i \leq M_i, i=1,2, \ldots, N\right\} 
$, $\mathbf{k} \in \mathcal{I}_{\mathbf{M}}$ means that $\mathbf{k}=\left(k_1, k_2, \ldots, k_N\right)$ fulfills  $1 \leq k_i \leq M_i$, $\forall i$, and  $\mathbf{x}=\left(x_1, x_2, \ldots, x_N\right) \in \Omega \subset \mathbb{R}^N$. 						
	Applying the operator \( A^{-2}BA^{\widetilde{\Gamma}} \) to \eqref{eq22}, we get
\begin{equation}\label{eq56}
	\begin{aligned}
		A^{-2} & BA^{\widetilde{\Gamma}} f(\mathbf{x}) = \sum_{\mathbf{k} \in \mathcal{I}_{\mathbf{M}}} \hat{f}(\mathbf{k}) (k_i\pi)^{2{\widetilde{\Gamma}}-3} \times \\ & \prod_{\substack{j=1 \\ j \neq i}}^N \sin \left( k_j \pi x_j \right)  \left( \cos(k_i\pi x_i) - 1 + a_1 x_i + \frac{k_i^2 \pi^2}{2} x_i^2 + \frac{k_i^2 \pi^2}{6} (\cos(k_i\pi) - 1) x_i^3 \right).	
	\end{aligned}					
\end{equation}	
	where
	\[
	a_1 = 1 - \cos(k_i\pi) - \frac{k_i^2 \pi^2}{2} - \frac{k_i^2 \pi^2}{6} (\cos(k_i \pi) - 1).
	\]
For the \( L^1(\Omega) \) norm, we apply the Riemann–Lebesgue lemma (see \cite{1684b3a8-48ea-3eb7-8a1e-11227e20bf15}) to obtain the following:
	\begin{align*}
		\left\| A^{-2}BA^{1-\zeta} f(\mathbf{x}) \right\|_1 &C \int_{\Omega}   \sum_{\mathbf{k} \in \mathcal{I}_{\mathbf{M}}}  \left| \hat{f}(\mathbf{k}) \right|  (k_i\pi)^{-2\zeta-1}  \mathrm{d}\mathbf{x} \leq  C \pi^{-2\zeta-1} \sum_{\mathbf{k} \in \mathcal{I}_{\mathbf{M}}} k_i^{-2\zeta-1}  .
	\end{align*}
Using \(  \hat{f}(\mathbf{k}) \in l^2(\Omega) \) (see \cite{garfken67:math}), we get:
	\begin{align*}
		\left\| A^{-2}BA^{\frac{5}{4}-\zeta} f(\mathbf{x}) \right\|_2^2 &\leq C \int_{\Omega}  \sum_{\mathbf{k} \in \mathcal{I}_{\mathbf{M}}}  \hat{f}(\mathbf{k})^2   \sum_{\mathbf{k} \in \mathcal{I}_{\mathbf{M}}}  (k_i\pi)^{-4\zeta-1}   \mathrm{d}\mathbf{x} \leq  C \sum_{\mathbf{k} \in \mathcal{I}_{\mathbf{M}}}  \hat{f}(\mathbf{k})^2   \sum_{\mathbf{k} \in \mathcal{I}_{\mathbf{M}}}  (k_i\pi)^{-4\zeta-1}  .
	\end{align*}
From the definition of the Fourier coefficients, it follows that \( \|f_k\|_{\infty} \leq C. \) Thus, an upper bound for \( A^{-2}BA^{1-\zeta} f(x) \) can be derived as follows:
	\begin{align*}
	\left\| A^{-2}BA^{1-\zeta} f(\mathbf{x}) \right\| \leq \sum_{\mathbf{k} \in \mathcal{I}_{\mathbf{M}}}  \left| \hat{f}(\mathbf{k}) \right|  (k_i\pi)^{-2\zeta-1}   \leq  C   \sum_{\mathbf{k} \in \mathcal{I}_{\mathbf{M}}}    (k_i\pi)^{-2\zeta-1}.
	\end{align*}
As \( M_1, \ldots, M_N \to \infty \), the desired outcomes for various norms are attained.
\end{proof}
 	\section{Numercal investigations}\label{numericalchap2}
In this section, we present numerical results to illustrate the error bounds established in Theorem \ref{theo4}.

\subsection{Implementation}
The implementation of exponential Runge--Kutta methods involves approximating the application of a matrix function to a vector. This is typically done by computing a weighted sum of functions \( \varphi_i(-\tau A) \)applied to vectors, written as \( \sum_{i=1}^k \varphi_i(-\tau A) v_k \).  To improve efficiency, an augmented matrix approach is used (see \cite{doi:10.1137/100788860}). This transforms the sum into a single matrix exponential acting on a vector, represented as \( \exp(-\tau \tilde{A}) V_0 \). 

There are a number of possible exponential Runge--Kutta methods that fulfill all the required order conditions. In this work, we choose a representative method, which is given by the following Butcher tableau:
\def\arraystretch{1.4}
$$
\begin{array}{c|cccc}
	0 & & \\
	\frac{1}{2} & \frac{1}{2} \varphi_{1,2} & \\
	\frac{1}{2} & \frac{1}{2}\varphi_{1,3} - \frac{1}{2}\varphi_{2,3} & \frac{1}{2} \varphi_{2,3} & \\
	1 & \varphi_{1,4} - 2\varphi_{2,4} &  \varphi_{1,4} & -\varphi_{1,4} + 2\varphi_{2,4} \\	
	\hline & \varphi_1-3 \varphi_2+4\varphi_3 & 0 & 4 \varphi_2-8 \varphi_3 & -\varphi_2 +  4 \varphi_3
\end{array} .
$$
The RK4 method with a sufficiently small time step is applied to find reference solution.

\subsection{Numerical results}
As in \cite{HAU2024}, the proof presented here also applies to any finite dimension. Therefore, it suffices to consider a one-dimensional domain, as extending the numerical experiments to higher dimensions is not particularly meaningful. We consider the following linear advection-diffusion problem 
\begin{equation}\label{problem2}
	\begin{aligned}
		& \partial_t  u(t,x) - 0.2 \partial_{xx} u (t,x)  =\partial_x u(t,x) ,\qquad (t,x)\in[0,1]\times[0,1], 
	\end{aligned}	
\end{equation} 
 subject to homogeneous Dirichlet boundary conditions. Since the operator \( A \) inherently incorporates these boundary conditions, the operators \( A \) and \( B \) do not commute. Additionally, the operators \( A^{-\frac{1}{2}}B \) and \( BA^{-\frac{1}{2}} \) are bounded and adhere to Assumption \ref{ass2} with \( \gamma = \frac{1}{2} \). For the discretization, we apply a standard second-order finite difference scheme to approximate the diffusion term \( -0.2 \partial_{xx} \) and the advection term \( \partial_{x} \), leading to the corresponding matrices \( A \) and \( B \). The computational domain is discretized using \( n=199 \) inner grid points. 

\begin{figure}
	\begin{center}
		\includegraphics[width=0.7\textwidth]{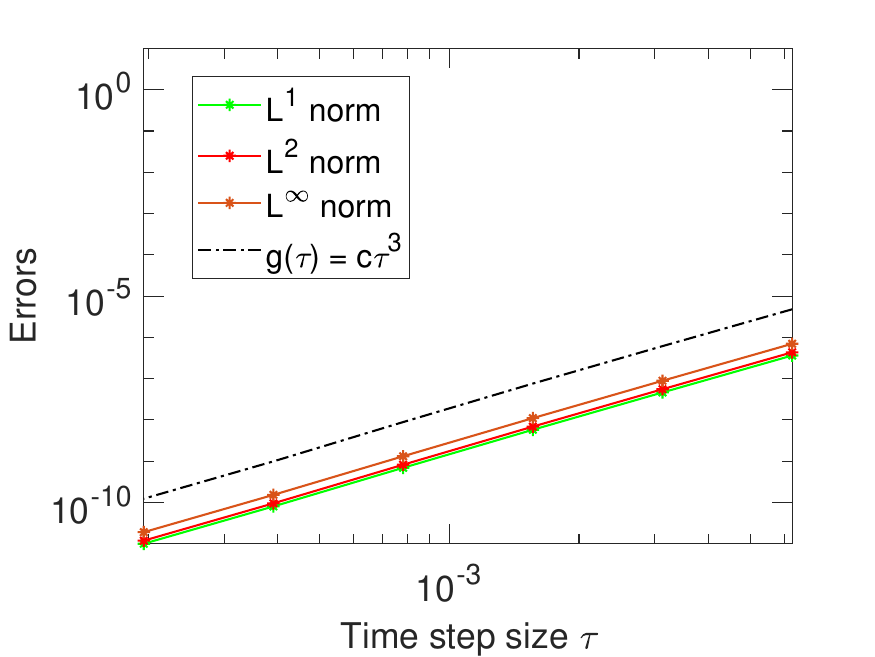}
		\caption{ The global error is shown as a function of the time step \( \tau \). }
		\label{fig1}
	\end{center}
\end{figure}

\begin{figure}
	\begin{center}
		\includegraphics[width=0.7\textwidth]{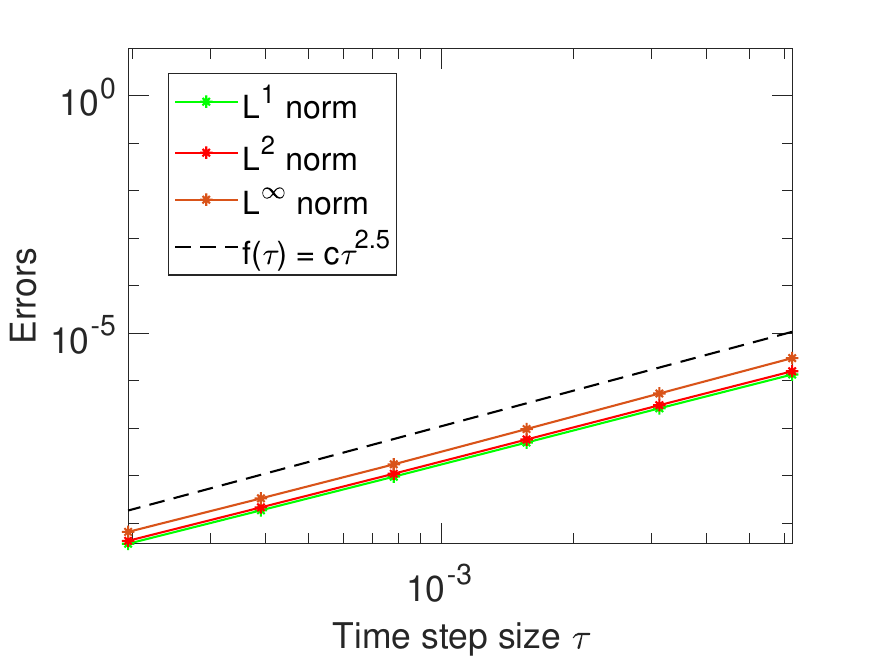}
		\caption{ The global error is shown as a function of the time step \( \tau \). }
		\label{fig2}
	\end{center}
\end{figure}

We choose the initial data \( 64x^3(1-x)^3 = u(0,x)   \in \mathcal{D}(A^2) \) to fulfill part of the requirements of Theorem \ref{theo4}. For the numerical simulations, we choose the abstract Banach spaces \( L^1(\Omega) \), \( L^2(\Omega) \), and \( L^{\infty}(\Omega) \) . The numerical results for various norms are presented in Figure \ref{fig1}, and they align perfectly with the Theorem \ref{theo4}. Additionally, we perform the ETD3RK method (see \cite{HO2005}) without avoiding order reduction, and the result is shown in Figure \ref{fig2}. As expected, the method exhibits order reduction from three to two and a half. This confirms the effectiveness of the approach in this paper to avoid order reduction.

\section{Conclusion}\label{concluchap2}
We have developed a fourth-stage Rung--Kutta method that satisfies all required order conditions under stronger assumptions on the initial data, applied to the initial value problem \eqref{4121} involving non-commuting  unbounded operators \( A \) and  \( B \). A thorough error analysis, conducted within the framework of semigroups in an abstract Banach space, demonstrates that the proposed method avoids order reduction. Higher-order methods could be further developed; however, this may involve more complicated modifications.


\bibliography{references}

\begin{thebibliography}{10}

\bibitem{doi:10.1137/100788860}
A.~H. Al-Mohy and N.~J. Higham.
\newblock Computing the action of the matrix exponential, with an application
  to exponential integrators.
\newblock {\em SIAM Journal on Scientific Computing}, 33(2):488--511, 2011.

\bibitem{garfken67:math}
G.~Arfken.
\newblock {\em Mathematical Methods for Physicists}.
\newblock Academic Press, {Inc.}, 1985.

\bibitem{1684b3a8-48ea-3eb7-8a1e-11227e20bf15}
S.~Bochner and K.~Chandrasekharan.
\newblock {\em Fourier Transforms. (AM-19)}.
\newblock Princeton University Press, 1949.

\bibitem{10.1007/s10543-013-0446-0}
M.~Caliari, P.~Kandolf, A.~Ostermann, and S.~Rainer.
\newblock Comparison of software for computing the action of the matrix
  exponential.
\newblock {\em BIT}, 54(1):113–128, 2014.

\bibitem{CROUSEILLES2020109688}
N.~Crouseilles, L.~Einkemmer, and J.~Massot.
\newblock Exponential methods for solving hyperbolic problems with application
  to collisionless kinetic equations.
\newblock {\em Journal of Computational Physics}, 420:109688, 2020.

\bibitem{DEKA2023101302}
P.~J. Deka, L.~Einkemmer, and M.~Tokman.
\newblock Le{XI}nt: Package for exponential integrators employing {L}eja
  interpolation.
\newblock {\em SoftwareX}, 21:101302, 2023.

\bibitem{henry1981geometric}
D.~Henry.
\newblock {\em Geometric Theory of Semilinear Parabolic Equations}.
\newblock Lecture notes in mathematics. Springer-Verlag, 1981.

\bibitem{HAU2024}
T.-H. Hoang.
\newblock Order reduction of exponential {R}unge--{K}utta methods:
  Non-commuting operators.
\newblock 2024.
\newblock Submitted to [Numerische Mathematik].

\bibitem{HO2005}
M.~Hochbruck and A.~Ostermann.
\newblock Explicit exponential {R}unge--{K}utta methods for semilinear
  parabolic problems.
\newblock {\em SIAM Journal on Numerical Analysis}, 43:1069--1090, 2005.

\bibitem{HOCHBRUCK2005323}
M.~Hochbruck and A.~Ostermann.
\newblock Exponential {R}unge–{K}utta methods for parabolic problems.
\newblock {\em Applied Numerical Mathematics}, 53(2):323--339, 2005.

\bibitem{HO2010}
M.~Hochbruck and A.~Ostermann.
\newblock Exponential integrators.
\newblock {\em Acta Numerica}, 19:209--286, 2010.

\bibitem{pazy1983semigroups}
A.~Pazy.
\newblock {\em Semigroups of Linear Operators and Applications to Partial
  Differential Equations}.
\newblock Applied Mathematical Sciences. Springer, 1983.

\end{thebibliography}

\end{document}